# Hypergeometric expansions of the general Heun function governed by two-term recurrence relations

T.A. Ishkhanyan[1,2] and A.M. Ishkhanyan[3,4]

[1]Moscow Institute of Physics and Technology, Dolgoprudny, 141700 Russia
[2]Institute for Physical Research, NAS of Armenia, Ashtarak, 0203 Armenia
[3]Russian-Armenian University, H. Emin 123, 0051, Yerevan, Armenia
[4]Institute of Physics and Technology, National Research Tomsk Polytechnic University, Tomsk, 634050 Russia

E-mail: aishkhanyan@gmail.com

We show that there exist infinitely many particular choices of parameters for which the three-term recurrence relations governing the expansions of the solutions of the general Heun equation in terms of the Gauss hypergeometric functions become two-term. In these cases the coefficients are explicitly written in terms of the gamma functions.

## 1. Introduction

The Heun equation [1] is the most general ordinary linear second-order Fuchsian differential equation having four regular singular points. The special functions emerging from this equation as well as from its four confluent forms [2-4] are widely faced in current fundamental and applied physics and mathematics (see, e.g., [2-4] and references therein). These functions generalize the hypergeometric, Lamé, Mathieu, spheroidal wave and many other known special functions, and for this reason they are believed to gradually become a part of the standard special functions of mathematical physics in the near future [3,4]. Currently, the Heun equations are encountered in so many branches of classical and non-classical physics ranging from non-Newtonian liquid mechanics, rheology, surface physics, polymer physics, atomic and nuclear physics to general gravity, astronomy and cosmology, e.g., the dislocation theory and mass-step problem in quantum mechanics, time-dependent few-state models in quantum optics, quantum two-center problem in molecular physics, surface polaritons, lattice systems in statistical mechanics, theory of black holes, etc., that it is difficult to give a classification of all relevant problems (see, e.g., [5-17]).

However, the Heun equations are much less studied than their predecessors, first of all, the hypergeometric relatives. A major reason for the slow progress in the development of the theory of the Heun equations is that the series solutions (either in terms of powers or in terms of the functions of hypergeometric class) in this case are governed by three-term recurrence relations between successive coefficients of expansions [2-4], instead of two-term ones appearing in the hypergeometric case [18-20]. For this reason, it was thought that the expansion coefficients are no longer determined explicitly.

However, it has recently been shown that there exist infinitely many particular choices of the involved parameters for which the recurrence relations for power-series expansions become two-term [21,22]. In these cases the solution of the Heun equation can be written either as a linear combination of a finite number of the hypergeometric functions or in terms of a single generalized hypergeometric function [21,22].

A further result we report here is that if the solutions of the general Heun equation are expanded in terms of the Gauss hypergeometric functions there also exist infinitely many particular choices of the involved parameters for which the governing three-term recurrence relations for expansion coefficients become two-term. In these cases the expansion coefficients are explicitly written in terms of the gamma functions.

## 2. Results

The general Heun equation written in its canonical form is [1]

$$\frac{d^2u}{dz^2} + \left(\frac{\gamma}{z} + \frac{\delta}{z-1} + \frac{\varepsilon}{z-a}\right)\frac{du}{dz} + \frac{\alpha\beta z - q}{z(z-1)(z-a)}u = 0, \qquad (1)$$

where the parameters satisfy the Fuchsian relation $1+\alpha+\beta = \gamma+\delta+\varepsilon$ in order to warrant that the singularity at infinity is regular. An expansion of the solution of this equation in terms of the Gauss hypergeometric functions is written as [15]

$$u = \sum_{n=0}^{\infty} c_n \cdot {}_2F_1(\alpha, \beta; \gamma + \varepsilon + n; z), \qquad (2)$$

where the expansion coefficients obey the three-term recurrence relation

$$R_n c_n + Q_{n-1} c_{n-1} + P_{n-2} c_{n-2} = 0 \qquad (3)$$

with
$$R_n = (1-a)n(\varepsilon + \gamma + n - 1), \qquad (4)$$

$$Q_n = -R_n + a(1+n-\delta)(n+\varepsilon) + (a\alpha\beta - q), \qquad (5)$$

$$P_n = -\frac{a}{n+\varepsilon+\gamma}(n+\varepsilon)(n+\varepsilon+\gamma-\alpha)(n+\varepsilon+\gamma-\beta). \qquad (6)$$

The expansion applies if $\alpha, \beta$ and $\gamma+\varepsilon$ are not zero or negative integers. The restrictions on $\alpha$ and $\beta$ assure that the hypergeometric functions are not polynomials of fixed degree.

Our result is that this recurrence admits two-term reductions for infinitely many particular choices of the involved parameters. These reductions are achieved by the following ansatz guessed from the results of [21,22]:

$$c_n = \left(\frac{1}{n}\frac{\prod_{k=1}^{N+2}(a_k-1+n)}{\prod_{k=1}^{N+1}(b_k-1+n)}\right)c_{n-1}, \tag{7}$$

where we put

$$a_1,...,a_N,a_{N+1},a_{N+2} = 1+e_1,...,1+e_N,\gamma+\varepsilon-\alpha,\gamma+\varepsilon-\beta, \tag{8}$$

$$b_1,...,b_N,b_{N+1} = e_1,...,e_N,\gamma+\varepsilon \tag{9}$$

with parameters $e_1,...,e_N$ to be defined later. Note that this ansatz implies that $e_1,...,e_N$ are not zero or negative integers.

The ratio $c_n/c_{n-1}$ is explicitly written as

$$\frac{c_n}{c_{n-1}} = \frac{(\gamma+\varepsilon-\alpha-1+n)(\gamma+\varepsilon-\alpha-1+n)}{(\gamma+\varepsilon-1+n)n}\prod_{k=1}^{N}\frac{e_k+n}{e_k-1+n}. \tag{10}$$

With this, the recurrence relation (3) is rewritten as

$$R_n\frac{(\gamma+\varepsilon-\alpha-1+n)(\gamma+\varepsilon-\beta-1+n)}{(\gamma+\varepsilon-1+n)n}\prod_{k=1}^{N}\frac{e_k+n}{e_k-1+n}+Q_{n-1}+$$
$$P_{n-2}\frac{(\gamma+\varepsilon-2+n)(n-1)}{(\gamma+\varepsilon-\alpha-2+n)(\gamma+\varepsilon-\beta-2+n)}\prod_{k=1}^{N}\frac{e_k-2+n}{e_k-1+n} = 0. \tag{11}$$

Substituting $R_n$ and $P_{n-2}$ from equations (4),(6) and cancelling the common denominator, this equation becomes

$$(1-a)(\gamma+\varepsilon-\alpha-1+n)(\gamma+\varepsilon-\beta-1+n)\prod_{k=1}^{N}(e_k+n)+$$
$$+Q_{n-1}\prod_{k=1}^{N}(e_k-1+n)-a(\varepsilon+n-2)(n-1)\prod_{k=1}^{N}(e_k-2+n) = 0. \tag{12}$$

This is a polynomial equation in $n$. Notably, it is of degree $N+1$, not $N+2$ as it may be guessed at first glance, because the highest-degree term proportional to $n^{N+2}$ identically vanishes. Hence, we have an equation of the form

$$\sum_{m=0}^{N+1} A_m(a,q;\alpha,\beta,\gamma,\delta,\varepsilon;e_1,...,e_N)n^m = 0. \tag{13}$$

Then, equating to zero the coefficients $A_m$ warrants the satisfaction of the three-term recurrence relation (3) for all $n$. We thus have $N+2$ equations $A_m = 0$, $m = 0, 1, .., N+1$, of which $N$ equations serve for determination of the parameters $e_{1,2,...,N}$ and the remaining two impose restrictions on the parameters of the Heun equation.

One of the latter restrictions is derived by calculating the coefficient $A_{N+1}$ of the term proportional to $n^{N+1}$ which is readily shown to be $2+N-\delta$. Hence,

$$\delta = 2 + N. \tag{14}$$

The second restriction imposed on the parameters of the Heun equation is checked to be a polynomial equation of the degree $N+1$ for the accessory parameter $q$.

Resolving the recurrence (10), the coefficients of expansion (2) are explicitly written in terms of the gamma functions as ( $c_0 = 1$ )

$$c_n = \frac{\Gamma(\gamma+\varepsilon)\Gamma(n+\gamma+\varepsilon-\alpha)\Gamma(n+\gamma+\varepsilon-\beta)}{n!\Gamma(\gamma+\varepsilon-\alpha)\Gamma(\gamma+\varepsilon-\beta)\Gamma(n+\gamma+\varepsilon)} \prod_{k=1}^{N} \frac{e_k+n}{e_k}, \quad n=0,1,2,.... \tag{15}$$

Since $\alpha, \beta$ and $e_1, ..., e_N$ are not zero or negative integers, with the help of the Fuchsian condition $1+\alpha+\beta = \gamma+\delta+\varepsilon$, it follows from this relation that $c_n$ may vanish only if $\alpha$ is a positive integer such that $0 < \alpha < 2+N$ or $\beta$ is a positive integer such that $0 < \beta < 2+N$. It can be checked that in these cases expansion (2) reduces to a rational function of $z$.

Here are the explicit solutions of the recurrence relation (3)-(6) for $N = 0, 1, 2$.

$N = 0$:
$$\delta = 2, \tag{16}$$

$$q = a\gamma + (\alpha-1)(\beta-1), \tag{17}$$

$$c_n = \frac{\Gamma(\gamma+\varepsilon)\Gamma(n+\gamma+\varepsilon-\alpha)\Gamma(n+\gamma+\varepsilon-\beta)}{n!\Gamma(\gamma+\varepsilon-\alpha)\Gamma(\gamma+\varepsilon-\beta)\Gamma(n+\gamma+\varepsilon)}. \tag{18}$$

$N = 1$:
$$\delta = 3, \tag{19}$$

$$q^2 + (\alpha-2)(\alpha-1)(\beta-2)(\beta-1) + a(4+2a-4\alpha-4\beta+3\alpha\beta)\gamma + \\ 2a^2\gamma^2 - q(4+a-3\alpha-3\beta+2\alpha\beta+3a\gamma) = 0, \tag{20}$$

$$e_1 = -q + a(1+\gamma) - 1 + (\alpha-1)(\beta-1), \tag{21}$$

$$c_n = \frac{\Gamma(\gamma+\varepsilon)\Gamma(n+\gamma+\varepsilon-\alpha)\Gamma(n+\gamma+\varepsilon-\beta)}{n!\Gamma(\gamma+\varepsilon-\alpha)\Gamma(\gamma+\varepsilon-\beta)\Gamma(n+\gamma+\varepsilon)} \frac{e_1+n}{e_1}. \tag{22}$$

$N = 2$:
$$\delta = 4, \tag{23}$$

$$q^3 - q^2(10 + 3\alpha(\beta - 2) - 6\beta + a(4 + 6\gamma)) +$$
$$q(33 - 2\alpha(2\beta - 5)(3\beta - 4) + \beta(11\beta - 40) + \alpha^2(11 + 3(\beta - 4)\beta) +$$
$$a^2(4 + \gamma(18 + 11\gamma)) + 2a(6 - 4\beta + 15\gamma - 11\beta\gamma + \alpha(2\beta + 6\beta\gamma - 11\gamma - 4))) -$$
$$(\alpha - 1)(\alpha - 2)(\alpha - 3)(\beta - 1)(\beta - 2)(\beta - 3) - \quad (24)$$
$$2a(18 + 6a^2 + 9(-3 + \alpha)\alpha - 27\beta + (37 - 11\alpha)\alpha\beta + (9 + \alpha(3\alpha - 11))\beta^2 +$$
$$a(9 - 9\alpha - 9\beta + 5\alpha\beta))\gamma + a^2(18\alpha + 18\beta - 18 - 18a - 11\alpha\beta)\gamma^2 - 6a^3\gamma^3 = 0,$$

$$e_1 + e_2 = -q + a(2 + \gamma) - 3 + (\alpha - 1)(\beta - 1), \quad (25)$$

$$\frac{(e_1 + 1)(e_2 + 1)}{e_1 e_2} = \frac{-q + a(\alpha - 3)(\beta - 3) + 3a\gamma}{(a - 1)(\alpha - 3)(\beta - 3)}, \quad (26)$$

$$c_n = \frac{\Gamma(\gamma + \varepsilon)\Gamma(n + \gamma + \varepsilon - \alpha)\Gamma(n + \gamma + \varepsilon - \beta)}{n!\,\Gamma(\gamma + \varepsilon - \alpha)\Gamma(\gamma + \varepsilon - \beta)\Gamma(n + \gamma + \varepsilon)} \frac{(e_1 + n)(e_2 + n)}{e_1 e_2}. \quad (28)$$

These results can readily be checked by direct verification of the recurrence (3)-(6). We conclude by noting that there exist expansions of the solutions of the Heun equation in terms of the Gauss hypergeometric functions the form of which differs from those applied in expansion (3) [23-25]. It can be checked that similar explicit solutions of the corresponding three-term recurrence relations can be constructed for these expansions as well.

Thus, we have shown that for the expansions of the solutions of the general Heun equation in terms of the Gauss hypergeometric functions there exist infinitely many particular choices of the involved parameters for which the three-term recurrence relations for expansion coefficients become two-term. The expansion coefficients are in these cases explicitly written in terms of the gamma functions.

**Acknowledgments**

This research has been conducted within the scope of the International Associated Laboratory IRMAS (CNRS-France & SCS-Armenia). The work has been supported by the Armenian State Committee of Science (SCS Grant No. 18RF-139), the Armenian National Science and Education Fund (ANSEF Grant No. PS-4986), and the project "Leading Russian Research Universities" (Grant No. FTI_24_2016 of the Tomsk Polytechnic University). T.A. Ishkhanyan acknowledges the support from SPIE through a 2017 Optics and Photonics Education Scholarship, and thanks the French Embassy in Armenia for a doctoral grant as well as the Agence universitaire de la Francophonie and Armenian State Committee of Science for a Scientific Mobility grant.